\documentclass[9pt,11pt]{article}%
\usepackage{amsfonts}
\usepackage[applemac]{inputenc}
\usepackage{amssymb}
\usepackage{makeidx}
\usepackage{graphicx}
\usepackage{amsmath}
\usepackage[a4paper]{geometry}
\usepackage[displaymath]{lineno}
\usepackage[singlespacing]{setspace}
\usepackage[normalem]{ulem}
\usepackage{lineno}
\usepackage{color}%
\providecommand{\U}[1]{\protect\rule{.1in}{.1in}}
\providecommand{\U}[1]{\protect\rule{.1in}{.1in}}

\newtheorem{theo}{Theorem}
\newtheorem{lem}[theo]{Lemma}
\newtheorem{prop}[theo]{Proposition}
\newtheorem{cor}[theo]{Corollary}

\newtheorem{exam}{Example}

\newenvironment{dem}[1][Proof]{\noindent \textbf{#1.} }{\ \rule{0.5em}{0.5em}}
\begin{document}

\title{A non-convex relaxed version\ of minimax theorems\thanks{Research supported by
MICIU of Spain and Universidad de Alicante (Contract Beatriz Galindo BEA-GAL
18/00205), Research Projects PGC2018-097960-B-C21 from MICINN of\ Spain and
AICO/2021/165 of Generalitat Valenciana.}}
\author{M.I.A.\ Ghitri\thanks{e-mail: moghitri@gmail.com}, A. Hantoute\thanks{e-mail:
hantoute@us.es (corresponding author)}\\{\scriptsize Department of Mathematics, University of\ Alicante, Spain}}
\date{}
\maketitle

\begin{abstract}
Given a subset $A\times B$ of a locally convex space $X\times Y$ (with $A$
compact) and a function $f:A\times B\rightarrow\overline{\mathbb{R}}$ such
that $f(\cdot,y),$ $y\in B,$ are concave and upper semicontinuous, the minimax
inequality $\max_{x\in A}\inf_{y\in B}f(x,y)\geq\inf_{y\in B}\sup_{x\in A_{0}%
}f(x,y)$ is shown to hold provided that $A_{0}$ be the set of $x\in A$ such
that $f(x,\cdot)$ is proper, convex and lower semi-contiuous. Moreover, if in
addition $A\times B\subset f^{-1}(\mathbb{R})$, then we can take as $A_{0}$
the set of $x\in A$ such that $f(x,\cdot)$ is convex. The relation to Moreau's
biconjugate representation theorem is discussed, and some applications
to\ convex duality are provided.



\end{abstract}

\section{Introduction\label{Sect1}}

Given a function $f:X\times Y\rightarrow\overline{\mathbb{R}},$ defined on the
Cartesian product of a locally convex space $X$\ and a linear space $Y,$ and
two convex sets $A\subset X,$ $B\subset Y,$ an enough general minimax theorem
(\cite{Fa53}) ensures that $\inf_{y\in B}\sup_{x\in A}f(x,y)=\max_{x\in A}%
\inf_{y\in B}f(x,y)$\ or, equivalently,
\begin{equation}
\inf_{y\in B}\sup_{x\in A}f(x,y)\leq\max_{x\in A}\inf_{y\in B}f(x,y).
\label{hh2}%
\end{equation}
Standard requirements for the validity of this minimax equality are the
compactness of the convex set $A,$ the concavity and upper semicontinuity
of\ the functions\ $f(\cdot,y)$ when $y$ ranges the whole set $B,$ together
with\ the following conditions,\
\begin{equation}
A\times B\subset f^{-1}(\mathbb{R}),\text{ and }f(x,\cdot)\text{ is convex for
all }x\in A. \label{hh}%
\end{equation}
First versions of the minimax theorem date back to \cite{VN59} and since then,
has not ceased to be improved. Several interesting generalizations of the
above theory have been proposed, either by relaxing the underlying linear
structure or convexity/concavity assumptions. We refer to \cite{Si95, Si58}
for an exhaustive presentation of many variants and extensions of this
theorem.\ In this paper, we are concerned with the following two questions
dealing with the possible relaxation of the two conditions in (\ref{hh}),
which to the best of our knowledge have not been considered before.

(1) Firstly, can the hypothesis on the finiteness of the function $f,$
$A\times B\subset f^{-1}(\mathbb{R}),$ be removed? Clearly, avoiding this
restrictive condition would widen the applicability of the theorem above\ to
large families of functions. For instance, to all convex/convex functions
taking infinite ($+\infty/-\infty$) values, which are of frequent use in
convex optimization and elsewhere.

(2) Does the inequality in (\ref{hh2}) remain true if the functions
$f(x,\cdot)$ are not\ convex for all $x\in A$ but, instead, only on a subset
of $A?$ In other words, does (\ref{hh2}) hold if\ $A$ is replaced with the set
$A_{0}:=\{x\in A:f(x,\cdot)$ convex$\}$?$\ $

To answer these two questions, we propose\ in this paper two relaxed variants
of the above minimax (inequality) theorem. For any given function $f:X\times
Y\rightarrow\overline{\mathbb{R}}$ such that $f(\cdot,y)$ are concave and
upper semicontinuous, for all $y\in B,$ and\ non-empty convex sets $A\subset
X,$ $B\subset Y,$ with $A$ being compact, we establish the following two results.

\begin{description}
\item[Theorem (see Theorem \ref{MM1})] If $A_{0}:=\{x\in A:f(x,\cdot)\in
\Gamma_{0}(Y)\},$ then
\[
\inf_{y\in Y}\sup_{x\in A_{0}}f(x,y)\leq\max_{x\in A}\inf_{y\in Y}f(x,y).
\]

\item[Theorem (see Theorem \ref{MMB})] If $A\times B\subset f^{-1}%
(\mathbb{R})$ and $A_{1}:=\{x\in A:f(x,\cdot)$ is convex$\},$ then
\[
\inf_{y\in B}\sup_{x\in A_{1}}f(x,y)\leq\max_{x\in A}\inf_{y\in B}f(x,y).
\]

\end{description}

The proofs of these two results are new, and are based on the Moreau
representation theorem of the biconjugate function (\cite{Mo03}), which is
itself a kind of a minimax theorem. The relationship between the minimax
theorem and the mentioned Moreau theorem had already been recognized
(\cite{Mo64, Ro64}), and in this paper we prove that they are indeed equivalent.

The above minimax theorems are of frequent use in optimization and convex
duality, we refer to \cite{CHL16, HL22, HMSV95, LiNg11, LV12} and references
therein for\ applications to subdifferential calculus of the supremum
functions. For the sake of motivation, we give the following\ example (see
Example \ref{exam} and Corollary \ref{cor} for the details).

\textbf{Example: }Given a finite family of lower semicontinuous convex
functions $\{\varphi_{1},$ $\ldots,$ $\varphi_{k}\}$ and $\varphi:=\max_{1\leq
i\leq k}\varphi_{i},$ we consider the optimization problem
\[
\inf_{x\in X}\varphi(x)\text{ }\left(  =\inf_{x\in X}\max_{\lambda\in
\Delta_{k}}\sum_{1\leq i\leq k}\lambda_{i}\varphi_{i}(x)\right)  ,\
\]
where $\Delta_{k}$ is the $k$-simplex. If all the $\varphi_{i}$'s are proper
(i.e., they do not take the value $-\infty$), then we can apply the standard
minimax theorem (\ref{hh2}) to the function $f(\lambda,x):=\sum_{1\leq i\leq
k}\lambda_{i}\varphi_{i}(x)$ and the sets $A:=\Delta_{k},$
$B=\operatorname*{dom}\varphi$ (the effective domain of $\varphi$) to
ensure\ that
\begin{equation}
\inf_{x\in X}\varphi(x)=\max_{\lambda\in\Delta_{k}}\inf_{x\in X}\text{{}}%
\sum_{1\leq i\leq k}\lambda_{i}\varphi_{i}(x) \label{df}%
\end{equation}
Otherwise, if the $\varphi_{i}$'s are not all proper, the function $f$ can
take the value $-\infty,$ and (\ref{hh2}) could not be applied, at least
directly.\ However, we can use the first\ theorem above (Theorem
\ref{MM1}),\ appealing to\ the set $A_{0}:=\{\lambda\in\Delta_{k}:\lambda
_{i}>0$ for all $i\}$ (observe that $0f_{i}$ is the indicator function of
$\operatorname*{dom}f_{i},$ and this function is not necessarily lower
semicontinuous). Then we obtain
\[
\inf_{x\in X}\max_{\lambda\in A_{0}}\sum_{1\leq i\leq k}\lambda_{i}\varphi
_{i}(x)\leq\max_{\lambda\in\Delta_{k}}\inf_{x\in X}\sum_{1\leq i\leq k}%
\lambda_{i}\varphi_{i}(x),
\]
and (\ref{df}) again follows as $\max_{\lambda\in A_{0}}\sum_{1\leq i\leq
k}\lambda_{i}\varphi_{i}(x)=\max_{\lambda\in\Delta_{k}}\sum_{1\leq i\leq
k}\lambda_{i}\varphi_{i}(x).$

Many developments have been made recently in the topic of minimax theorems,
which could also be considered in our case for further developments. We refer,
for instance, to \cite{BS18} for extensions to abstract convex functions, to
\cite{Ga14} for relaxed convexity conditions using partial data, to
\cite{JLL09, Ro22} for application to alternative theorems, to \cite{Ri17} and
\cite{Sa18} for new topological minimax theorems. Other extensions are
given\ in \cite{La17}, \cite{Mo12} and \cite{Su18} among other achievements.

The paper is organized as follows. In section \ref{Sect2} we introduce the
main notations and results that are needed in the sequel. Section \ref{top}
reviews some topological minimax theorems. The extended real-valued version of
the minimax theorem, Theorem \ref{MM1}, is given in section \ref{minimax1},
whereas the second variant given in Theorem \ref{MMB}, is given in section
\ref{minimax2}.

\section{Preliminaries\label{Sect2}}

Throughout the paper, $X$ stands for a real (separated) locally convex spaces
(lcs, for short). The pairing between $X$ and its\ topological dual $X^{\ast}$
is\ denoted\ by $(x^{\ast},x)\in X^{\ast}\times X\mapsto\langle x^{\ast
},x\rangle:=\left\langle x,x^{\ast}\right\rangle :=x^{\ast}(x)$. The zero
vector is denoted by $\theta$, and the family of closed convex balanced
neighborhoods of $\theta,$\textbf{ }called\textbf{ }$\theta$-neighborhoods, is
denoted $\mathcal{N}_{X}.$ We use the notation $\overline{\mathbb{R}%
}:=\mathbb{R}\cup\{-\infty,+\infty\}$ and $\mathbb{R}_{\infty}:=\mathbb{R}%
\cup\{+\infty\}$, and adopt the conventions\emph{ }$\left(  +\infty\right)
+(-\infty)=\left(  -\infty\right)  +(+\infty)=+\infty,$ $0.(+\infty)=+\infty,$
and $0(-\infty)=0.$

Given a set $A\subset X,$ by $\operatorname*{co}A$ we denote the convex
hull\emph{ }of the set $A$, while $\operatorname*{cl}A$ (and sometimes
$\overline{A})$ is used for denoting the closure of $A$; in particular,
$\overline{\operatorname*{co}}A:=\operatorname*{cl}(\operatorname*{co}A)$. The
polar set of $A$ is $A^{\circ}:=\{x^{\ast}\in X^{\ast}:\left\langle x^{\ast
},x\right\rangle \leq1$ for all $x\in A\}.$ The indicator function of $A$ is
the function $\mathrm{I}_{A}:X\rightarrow\mathbb{R}_{\infty}$ defined by
$\mathrm{I}_{A}(x):=0$ if $x\in A,$ and $\mathrm{I}_{A}(x):=+\infty$ otherwise.

Given a function $f:X\longrightarrow\overline{\mathbb{R}}$, by
$\operatorname*{dom}f:=\{x\in X:\ f(x)<+\infty\}$ and $\operatorname*{epi}%
f:=\{(x,\lambda)\in X\times\mathbb{R}:\ f(x)\leq\lambda\},$ we denote the
(effective) domain and epigraph of $f$, respectively. We say that $f$ is
proper if $\operatorname*{dom}f\neq\emptyset$ and $f(x)>-\infty$ for all $x\in
X,$ lower semicontinuous (lsc, for short) if $\operatorname*{epi}f$ is closed,
and convex if $\operatorname*{epi}f$ is convex.\ The closed and the closed
convex hulls of $f$ are, respectively,\ the functions $\operatorname*{cl}f$
(or $\bar{f})$ and $\overline{\operatorname*{co}}f$
satisfying\ $\operatorname*{epi}\bar{f}=\operatorname*{cl}(\operatorname*{epi}%
f)$ and $\operatorname*{epi}(\overline{\operatorname*{co}}f)=\overline
{\operatorname*{co}}(\operatorname*{epi}f).$ We shall denote by $\Gamma
_{0}(X)$ the family of proper, convex and lsc functions defined on $X.$ It is
known that a lsc convex function which takes somewhere the value $-\infty,$ it
can only take infinite values $(+\infty$ and $-\infty).$

The conjugate of $f$ is the lsc convex function $f^{\ast}:X^{\ast}%
\rightarrow\overline{\mathbb{R}}\ $defined by
\[
f^{\ast}(x^{\ast}):=\sup\{\left\langle x^{\ast},x\right\rangle -f(x):\ x\in
X\}.
\]
The function $f^{\ast}$ is proper if and only if $f$ is, provided that $f$ is
lsc and convex. Moreover, due to the relations
\begin{equation}
\inf f=\inf(\operatorname*{cl}f)=\inf(\operatorname*{co}f)=\inf(\overline
{\operatorname*{co}}f), \label{ehuls}%
\end{equation}
the function $f^{\ast}$ does not distinguish between the function and its
convex hulls;\ that is,
\begin{equation}
f^{\ast}=(\operatorname*{co}f)^{\ast}=(\overline{\operatorname*{co}}f)^{\ast}.
\label{tricy}%
\end{equation}
Moreover, due to Moreau's theorem (\cite{Za02}), provided that $\overline
{\operatorname*{co}}f$ is proper we have
\begin{equation}
f^{\ast\ast}=\overline{\operatorname*{co}}f, \label{moreau}%
\end{equation}
where $f^{\ast\ast}:X\rightarrow\overline{\mathbb{R}}$ stands for (the
restriction to $X$ of) the conjugate of $f^{\ast}.$ The conjugation operation
obeys many nice calculus rules; for instance, given a family of functions
$f_{t}:X\rightarrow\overline{\mathbb{R}},$ $t\in T,$ the
following\ relation\ always holds
\begin{equation}
\left(  \inf\limits_{t\in T}f_{t}\right)  ^{\ast}=\sup\limits_{t\in T}%
f_{t}^{\ast}, \label{conv3}%
\end{equation}
whereas (\ref{moreau}) gives rise to
\begin{equation}
\left(  \sup\limits_{t\in T}f_{t}\right)  ^{\ast}=\overline{\operatorname*{co}%
}\left(  \inf\limits_{t\in T}f_{t}^{\ast}\right)  ; \label{infsupcon}%
\end{equation}
the latter\ being true when $f_{t}$ and $\sup_{t\in T}f_{t}\in\Gamma_{0}(X).$
A related concept is the $\varepsilon$\emph{-}subdifferential\emph{ }of $f,$
$\varepsilon\in\mathbb{R}$, given by
\[
\partial_{\varepsilon}f(x):=\{x^{\ast}\in X^{\ast}:f(y)\geq f(x)+\langle
x^{\ast},y-x\rangle-\varepsilon,\text{\ for all }y\in X\},
\]
with $\partial_{\varepsilon}f(x):=\emptyset$ if $x\notin\operatorname*{dom}f$
or $\varepsilon<0.$\ The set $\partial f(x):=\partial_{0}f(x)$ is the
subdifferential of $f$ at $x.$

We recall the following known fact that will be used later (see \cite{HLZ08}):
If\ $\{f_{t}:t\in T\}$ is a non-empty family of convex functions, and
$f:=\sup_{t\in T}f_{t}$, are such that $\operatorname*{aff}%
(\operatorname*{dom}f_{t})=\operatorname*{aff}(\operatorname*{dom}f),$ for all
$t\in T,$ and $f_{\mid\operatorname*{aff}(\operatorname*{dom}f)}$ is finite
and continuous somewhere in $\operatorname*{ri}(\operatorname*{dom}f),$ then
we have
\begin{equation}
\operatorname*{cl}f=\sup_{t\in T}\operatorname*{cl}f_{t}. \label{amin}%
\end{equation}

\section{A topological minimax theorem\label{top}}

The present section provides topologically flavored minimax-type results\ that
illustrate the essential arguments underlying the proof of minimax theorems,
namely, the utility of compactness and monotonicity-like properties with
respect to one of the variables. The results of this section will also be used later.

The first\ proposition\ analyzes\ the interchange between the lower limit and
the minimum up to some appropriate lsc regularizations. This result can be
obtained from \cite[Propostition 7.29]{RW98} but, for the sake of
completeness, a short proof is given here with\ a slight improvement\ of the
underlying upper semi-continuity assumption. Throughout this section, $X$
stands for a\ topological space.

\begin{prop}
\label{mmmonotonemain} Let $(\varphi_{i})_{i\in I}$ be a net of extended
real-valued functions defined on a topological space $X.$ For any nonempty
compact set $A\subset X,$\ we have\
\begin{equation}
\min_{x\in A}\text{ }\sup_{i\in I}\text{ }\operatorname*{cl}\left(
\inf_{i\preccurlyeq j}\varphi_{j}+\mathrm{I}_{A}\right)  (x)=\sup_{i\in
I}\text{ }\inf_{x\in A,\text{ }i\preccurlyeq j}\text{ }\varphi_{j}(x).
\label{ineqmmain}%
\end{equation}

\end{prop}

\begin{dem}
We denote $\psi_{i}:=\operatorname*{cl}\left(  \inf_{i\preccurlyeq j}%
\varphi_{j}+\mathrm{I}_{A}\right)  ,$ $i\in I,$ so that each $\psi_{i}$ is
lsc, obviously,\ and\ $(\psi_{i})_{i\in I}$ forms a non-decreasing net. Then,
taking into account the compactness of $A,$\ we pick\ a\ net $(x_{i}%
)_{i}\subset A$ such that $\psi_{i}(x_{i})=\min_{x\in A}\psi_{i}(x),$ for each
$i\in I,$ and $x_{i}\longrightarrow\bar{x}\in A$ without loss of generality
(w.l.o.g.~for short) on $I$. Thus, for each given $j\in I,$\ we have $\psi
_{j}(\bar{x})\leq\liminf_{i\in I}\psi_{j}(x_{i})$\ and, so,
\[
\psi_{j}(\bar{x})\leq\liminf_{i\in I}\psi_{j}(x_{i})\leq\liminf_{i\in I}%
\psi_{i}(x_{i})=\liminf_{i\in I}\min_{x\in A}\psi_{i}(x).
\]
Consequently, since the last inequality holds for all $j\in I,$\
\begin{align*}
\min_{x\in A}\sup_{i\in I}\psi_{i}(x)  &  =\min_{x\in A}\lim_{i\in I}\psi
_{i}(x)\leq\lim_{i\in I}\psi_{i}(\bar{x})=\sup_{i\in I}\min_{x\in A}\psi
_{i}(x)\\
&  \leq\sup_{i\in I}\inf_{x\in A}\left(  \inf_{i\preccurlyeq j}\varphi
_{j}+\mathrm{I}_{A}\right)  (x)=\sup_{i\in I}\inf_{x\in A,\text{
}i\preccurlyeq j}\varphi_{j}(x),
\end{align*}
and relation\ (\ref{ineqmmain}) follows\ as%
\begin{align*}
\sup_{i\in I}\inf_{x\in A,\text{ }i\preccurlyeq j}\varphi_{j}(x)  &
=\sup_{i\in I}\inf_{x\in X}\left(  \inf_{i\preccurlyeq j}\varphi
_{j}+\mathrm{I}_{A}\right)  (x)\\
&  =\sup_{i\in I}\inf_{x\in X}\operatorname*{cl}\left(  \inf_{i\preccurlyeq
j}\varphi_{j}+\mathrm{I}_{A}\right)  (x)=\sup_{i\in I}\inf_{x\in X}\psi
_{i}(x)\leq\min_{x\in A}\sup_{i\in I}\psi_{i}(x).
\end{align*}

\end{dem}

Next, we discusses some situations in which (\ref{ineqmmain}) simplifies
to\ the minimax relation $\sup_{i\in I}\min_{x\in A}\varphi_{i}(x)=\min_{x\in
A}\sup_{i\in I}\varphi_{i}(x);$ the validity of the latter in general is
discarded (see, e.g.,\ \cite[page 239]{RW98}). A first consequence of
Proposition \ref{mmmonotonemain} concerns non-decreasing nets $(\varphi
_{i})_{i\in I};$ that is,
\[
i_{1}\preccurlyeq i_{2}\Rightarrow\varphi_{i_{1}}\leq\varphi_{i_{2}}\text{
\ for all }i_{1},\text{ }i_{2}\in I,
\]
where \textquotedblleft$\preccurlyeq$\textquotedblright\ stands for\ the
partial order relation directing\ the index set $I.\ $

\begin{cor}
\label{mmmonotone} Assume that the net $(\varphi_{i})_{i\in I}$ in Proposition
\ref{mmmonotonemain} is non-decreasing.$\ $Then we have that
\begin{equation}
\min_{x\in A}\sup_{i\in I}\operatorname*{cl}(\varphi_{i}+\mathrm{I}%
_{A})(x)=\sup_{i\in I}\inf_{x\in A}\varphi_{i}(x), \label{ineqmm}%
\end{equation}
and consequently, provided that the $\varphi_{i}$'s are lsc,%
\[
\min_{x\in A}\sup_{i\in I}\varphi_{i}(x)=\sup_{i\in I}\min_{x\in A}\varphi
_{i}(x).
\]

\end{cor}

\begin{dem}
Under the current assumption, each of the functions\ $\inf_{i\preccurlyeq
j}\varphi_{j}$ $(=\varphi_{j}),$ $i\in I,$ is lsc and (\ref{ineqmmain})
straightforwardly gives rise to (\ref{ineqmm}).
\end{dem}

Convexity would make it possible to avoid the lower semi-continuity assumption
used in the previous corollary; as\ for the property of monotonicity, we
must\ await the following sections.

\begin{cor}
\label{mmmonotones} Let\ $(\varphi_{i})_{i\in I}$ be a non-decreasing net of
extended real-valued proper convex functions defined on a lcs space $X,$ and
let $A\subset X$ be a nonempty compact convex set. We assume that the function
$\varphi:=\sup_{i\in I}\varphi_{i}$ is finite and continuous at some point in
$\operatorname*{int}(A),$ assumed nonempty. Then we have that
\[
\inf_{x\in A}\sup_{i\in I}\varphi_{i}(x)=\sup_{i\in I}\inf_{x\in A}\varphi
_{i}(x),
\]

\end{cor}

\begin{dem}
The current convexity and continuity/interiority assumptions guarantee,
according to \cite[Corollary 9 and Lemma 15]{HLZ08}, that $\sup_{i\in
I}(\operatorname*{cl}(\varphi_{i}+\mathrm{I}_{A}))=\operatorname*{cl}\left(
\sup_{i\in I}\left(  \varphi_{i}+\mathrm{I}_{A}\right)  \right)  $ and
$\operatorname*{cl}(\varphi_{i}+\mathrm{I}_{A})=(\operatorname*{cl}\varphi
_{i})+\mathrm{I}_{A},$ for all $i\in I.\ $So,
\begin{align*}
\min_{x\in A}\sup_{i\in I}\operatorname*{cl}\left(  \varphi_{i}+\mathrm{I}%
_{A}\right)  (x)  &  =\min_{x\in X}\sup_{i\in I}\left(  \operatorname*{cl}%
(\varphi_{i}+\mathrm{I}_{A})\right)  (x)\\
&  =\min_{x\in X}\operatorname*{cl}\left(  \sup_{i\in I}\left(  \varphi
_{i}+\mathrm{I}_{A}\right)  \right)  (x)\\
&  =\inf_{x\in X}\sup_{i\in I}\left(  \varphi_{i}+\mathrm{I}_{A}\right)
(x)=\inf_{x\in A}\sup_{i\in I}\varphi_{i}(x),
\end{align*}
and the desired conclusion follows from (\ref{ineqmm}).
\end{dem}

One possible way to overcome the compactness assumption used above is to
compactify the given set $A$ and, consequently, use appropriate lsc
regularization of the $\varphi_{i}$'s. For this purpose, given a topological
space $X,$ we consider the Stone-\v{C}ech compactification $\beta(X)$
of\ $X$\ defined\ as the closure of $X$ in the (Hausdorff)\ compact product
topological space
\[
T:=\left[  0,1\right]  ^{\mathcal{C}(X,\left[  0,1\right]  )}\equiv
\{\gamma:\mathcal{C}(X,\left[  0,1\right]  )\rightarrow\left[  0,1\right]
\},
\]
where $\mathcal{C}(X,\left[  0,1\right]  )$ denotes\ the set of continuous
functions from $X$ to $\left[  0,1\right]  \ $(we refer to \cite{Mu00} for
more details on this compactification process). Indeed, $X$ is identified
with\ a subset of\ $T\mathbb{\ }$thanks to\ the mapping\ $\mathfrak{w}%
:X\rightarrow T$ defined as $\mathfrak{w}(x)\equiv\gamma_{x},$ where
\begin{equation}
\gamma_{x}(\varphi):=\varphi(x)\ \text{\ \ for all }\varphi\in\mathcal{C}%
(X,\left[  0,1\right]  ), \label{eval}%
\end{equation}
allowing us to\ set $X\equiv\mathfrak{w}(X).$ Recall that the convergence of a
given\ net\ $(\gamma_{j})_{j}$ to $\gamma$ in $T,$ written\ $\gamma
_{j}\rightarrow\gamma$, means that
\begin{equation}
\gamma_{j}(\varphi)\rightarrow\gamma(\varphi)\text{ \ \ for all }\varphi
\in\mathcal{C}(X,\left[  0,1\right]  ). \label{conve}%
\end{equation}
When\ $X$ is a Tychonoff space (i.e., completely regular and Hausdorff
(\cite{Mu00})), $\beta(X)$ is a Hausdorff space and the mapping $\mathfrak{w}$
is a homeomorphism between $X$ and $\mathfrak{w}(X);$ that is,
\begin{equation}
\gamma_{x_{j}}\rightarrow\gamma_{x}\text{ if and only if }x_{j}\rightarrow
x\text{ in }X, \label{hemeo}%
\end{equation}
for every $x\in X$ and every net $(x_{j})_{j}\subset X.$

\begin{cor}
\label{corcompact} Let\ $(\varphi_{i})_{i\in I}$ be a non-decreasing net of
extended real-valued functions defined on a topological\ space $X.$ Then we
have that
\[
\min_{\gamma\in\beta(X)}\sup_{i\in I}(\operatorname*{cl}\varphi_{i}%
)(\gamma)=\sup_{i\in I}\inf_{x\in X}\varphi_{i}(x),
\]
where $\beta(X)$ stands for\ the Stone-\v{C}ech compactification of $X$ and
$\operatorname*{cl}\varphi_{i}:\beta(X)\rightarrow\overline{\mathbb{R}},$
$i\in I,$ are defined by\
\[
(\operatorname*{cl}\varphi_{i})(\gamma):=\liminf_{\gamma_{x_{j}}%
\rightarrow\gamma,\text{ }(x_{j})_{j}\subset X}\varphi_{i}(x_{j}).
\]

\end{cor}

\begin{dem}
We apply Corollary\ \ref{mmmonotone} to the non-decreasing net $(\psi
_{i})_{i\in I}$ of extended real-valued functions defined on the topological
space $\beta(X)$ as
\[
\psi_{i}(\gamma):=\left\{
\begin{array}
[c]{ll}%
\varphi_{i}(x), & \text{if }\gamma=\gamma_{x},\text{ }x\in X,\\
+\infty, & \text{otherwise; }%
\end{array}
\right.
\]
hence, the closure of each $\psi_{i}$ with respect to the new topology on
$\beta(X)$ satisfies
\[
(\operatorname*{cl}\psi_{i})(\gamma)=\liminf_{\gamma_{x_{j}}\rightarrow
\gamma,\text{ }(x_{j})_{j}\subset X}\varphi_{i}(x_{j})=(\operatorname*{cl}%
\varphi_{i})(\gamma).
\]
Consequently,\ since $\beta(X)$ is compact by construction,
Corollary\ \ref{mmmonotone} yields
\[
\min_{\gamma\in\beta(X)}\sup_{i\in I}(\operatorname*{cl}\varphi_{i}%
)(\gamma)=\min_{\gamma\in\beta(X)}\sup_{i\in I}(\operatorname*{cl}\psi
_{i})(x)=\sup_{i\in I}\inf_{\gamma\in\beta(X)}\psi_{i}(\gamma)=\sup_{i\in
I}\inf_{x\in X}\varphi_{i}(x).
\]

\end{dem}

Corollary \ref{corcompact} is illustrated in the following result, presented
as an example, which will be used later when\ deriving the Moreau biconjugate
representation theorem from the minimax theorem. Note that although we have
used here the process of compactification, the final minimax
equality\ is\ expressed only using the original space $X.$

\begin{exam}
\label{exam} Let $X$ be a lcs space, $f\in\Gamma_{0}(X),$ and $x_{0}\in X.$
Then we have (the minimax\ inequality)
\begin{align}
\sup_{U\in\mathcal{N}_{X}}\inf_{x\in X}\left(  \sup_{x^{\ast}\in U^{\circ}%
}\{\left\langle x^{\ast},x_{0}-x\right\rangle +f(x)\}\right)   &  =\min_{x\in
X}\sup_{U\in\mathcal{N}_{X}}\left(  \sup_{x^{\ast}\in U^{\circ}}\{\left\langle
x^{\ast},x_{0}-x\right\rangle +f(x)\}\right) \label{bachir}\\
\text{ }(  &  =f(x_{0})).\nonumber
\end{align}

\end{exam}

\begin{dem}
We introduce\ the functions $\varphi_{U}:X\rightarrow\mathbb{R}_{\infty},$
$U\in\mathcal{N}_{X},$ given\ by
\[
\varphi_{U}(x):=\sup_{x^{\ast}\in U^{\circ}}\{\left\langle x^{\ast}%
,x_{0}-x\right\rangle +f(x)\},
\]
and endow the family $\mathcal{N}_{X}$ of $\theta$-neighborhoods in $X$ with
the partial order given by descending inclusions; that is,
\[
U_{1}\preccurlyeq U_{2},\text{ }U_{1},U_{2}\in\mathcal{N}_{X}\Leftrightarrow
U_{2}\subset U_{1},
\]
so that the net $(\varphi_{U})_{U\in\mathcal{N}_{X}}$ is non-decreasing.
Moreover, since each\ set $U^{\circ}$ is $w^{\ast}$-compact by\ Dieudonné's
Theorem (see, e.g., \cite[Theorem 1.1.8]{Za02}), we have that\ $\emptyset
\neq\operatorname*{dom}f\subset\operatorname*{dom}\varphi_{U}$ and, so,
$\varphi_{U}\in\Gamma_{0}(X).$\ Therefore, by Corollary \ref{corcompact},
\begin{equation}
\min_{\gamma\in\beta(X)}\sup_{U\in\mathcal{N}_{X}}(\operatorname*{cl}%
\varphi_{U})(\gamma)=\sup_{U\in\mathcal{N}_{X}}\inf_{x\in X}\varphi_{U}(x),
\label{hamdo}%
\end{equation}
where $\beta(X)$ stands for\ the Stone-\v{C}ech compactification of the lcs
$X$ and the function $\operatorname*{cl}\varphi_{U}:\beta(X)\rightarrow
\overline{\mathbb{R}},$ $U\in\mathcal{N}_{X},$ is defined as\
\[
(\operatorname*{cl}\varphi_{U})(\gamma):=\liminf_{\gamma_{x_{i}}%
\rightarrow\gamma,\text{ }(x_{i})_{i}\subset X}\varphi_{U}(x_{i}).
\]
In particular, since $X$ is completely regular (as $X$ is an (Hausdorff) lcs),
$\beta(X)$ is Hausdorff and the lower semicontinuity of each $\varphi_{U}$
yields, for every $x\in X,$
\begin{equation}
(\operatorname*{cl}\varphi_{U})(\gamma_{x})=\liminf_{\gamma_{x_{i}}%
\rightarrow\gamma_{x},\text{ }(x_{i})_{i}\subset X}\varphi_{U}(x_{i}%
)=\varphi_{U}(x). \label{karim}%
\end{equation}
Let us show that
\[
\min_{\gamma\in\beta(X)}\sup_{U\in\mathcal{N}_{X}}(\operatorname*{cl}%
\varphi_{U})(\gamma)=\sup_{U\in\mathcal{N}_{X}}\varphi_{U}(x_{0})=f(x_{0}).
\]
Indeed, if\ $\gamma\in\beta(X)\setminus\{\gamma_{x_{0}}\},$\ then
(\ref{hemeo}) together with\ the fact that $\beta(X)$ is Hausdorff gives rise
to\ some $V\in\mathcal{N}_{X}$ such that $x_{0}-x_{i}\in X\setminus V$
frequently for $i,$\ for all nets $(x_{i})_{i}\subset X$ such that
$\gamma_{x_{i}}\rightarrow\gamma.$ Hence, by the bipolar theorem
(\cite[Theorem 1.1.9]{Za02}), there exists some $x_{0}^{\ast}\in V^{\circ}%
$\ such that $\left\langle x_{0}^{\ast},x_{0}-x_{i}\right\rangle >1$
frequently for $i$ (for all nets $(x_{i})_{i}\subset X$ such that
$\gamma_{x_{i}}\rightarrow\gamma).$ Thus,\ using the lower semicontinuity of
the function $x\mapsto\sup_{x^{\ast}\in U^{\circ}}\{\left\langle x^{\ast
},x_{0}-x\right\rangle +f(x)\},$ we obtain
\begin{align*}
\sup_{U\in\mathcal{N}_{X}}(\operatorname*{cl}\varphi_{U})(\gamma)  &
=\sup_{U\in\mathcal{N}_{X}}\text{ }\liminf_{\gamma_{x_{i}}\rightarrow
\gamma,\text{ }(x_{i})_{i}\subset X}\text{ }\sup_{x^{\ast}\in U^{\circ}%
}\{\left\langle x^{\ast},x_{0}-x_{i}\right\rangle +f(x_{i})\}\\
&  \geq\sup_{U\in\mathcal{N}_{X}}\text{ }\sup_{x^{\ast}\in U^{\circ}}\text{
}\liminf_{\gamma_{x_{i}}\rightarrow\gamma,\text{ }(x_{i})_{i}\subset X}\text{
}\left(  \left\langle x^{\ast},x_{0}-x_{i}\right\rangle +f(x_{i})\right)  ,
\end{align*}
At the same time, since\ $f\in\Gamma_{0}(X),$ there are $a_{0}^{\ast}\in
X^{\ast}$ and $\alpha_{0}\in\mathbb{R}$ \ such that $f\geq\left\langle
a_{0}^{\ast},\cdot\right\rangle +\alpha_{0}$ (\cite[Theorem 2.2.6]{Za02}), and
the last inequality above simplifies to%
\[
\sup_{U\in\mathcal{N}_{X}}(\operatorname*{cl}\varphi_{U})(\gamma)\geq
\sup_{U\in\mathcal{N}_{X}}\text{ }\sup_{x^{\ast}\in U^{\circ}}\text{ }%
\liminf_{\gamma_{x_{i}}\rightarrow\gamma,\text{ }(x_{i})_{i}\subset
X}\text{\ }\left\langle x^{\ast}-a_{0}^{\ast},x_{0}-x_{i}\right\rangle
+\left(  \left\langle a_{0}^{\ast},x_{0}\right\rangle +\alpha_{0}\right)  .
\]
Now, given any\ $m\geq1,$ we have $m^{-1}V\in\mathcal{N}_{X}$ and
$mx_{0}^{\ast}\in mV^{\circ}=(m^{-1}V)^{\circ}.$ So, choosing\ $U$ small
enough such that $a_{0}^{\ast}+mV^{\circ}\subset U^{\circ},$ we obtain
$a_{0}^{\ast}+mx_{0}^{\ast}\in U^{\circ}$ and the last inequality yields
\begin{align*}
\sup_{U\in\mathcal{N}_{X}}(\operatorname*{cl}\varphi_{U})(\gamma)  &  \geq
\sup_{x^{\ast}\in a_{0}^{\ast}+mV^{\circ}}\text{ }\liminf_{\gamma_{x_{i}%
}\rightarrow\gamma,\text{ }(x_{i})_{i}\subset X}\text{ }\left\langle x^{\ast
}-a_{0}^{\ast},x_{0}-x_{i}\right\rangle +\left(  \left\langle a_{0}^{\ast
},x_{0}\right\rangle +\alpha_{0}\right) \\
&  \geq\liminf_{\gamma_{x_{i}}\rightarrow\gamma,\text{ }(x_{i})_{i}\subset
X}\text{ }\left\langle mx_{0}^{\ast},x_{0}-x_{i}\right\rangle +\left(
\left\langle a_{0}^{\ast},x_{0}\right\rangle +\alpha_{0}\right) \\
&  \geq m+\left(  \left\langle a_{0}^{\ast},x_{0}\right\rangle +\alpha
_{0}\right)  .
\end{align*}
In other words, by the arbitrariness of $m\geq1,$ we have that $\sup
_{U\in\mathcal{N}_{X}}(\operatorname*{cl}\varphi_{U})(\gamma)=+\infty,$\ and
the infimum in $\inf_{\gamma\in\beta(X)}\sup_{U\in\mathcal{N}_{X}%
}(\operatorname*{cl}\varphi_{U})(\gamma)$ is\ attained at $\gamma
=\gamma_{x_{0}}.$ Consequently, by combining (\ref{hamdo}) and (\ref{karim})
we infer that%
\begin{align*}
\sup_{U\in\mathcal{N}_{X}}\inf_{x\in X}\varphi_{U}(x)  &  =\min_{\gamma
\in\beta(X)}\sup_{U\in\mathcal{N}_{X}}(\operatorname*{cl}\varphi_{U}%
)(\gamma)=\min_{\gamma_{x}\in\beta(X)}\sup_{U\in\mathcal{N}_{X}}%
(\operatorname*{cl}\varphi_{U})(\gamma_{x})\\
&  =\min_{x\in X}\sup_{U\in\mathcal{N}_{X}}\varphi_{U}(x)=\sup_{U\in
\mathcal{N}_{X}}\varphi_{U}(x_{0})=f(x_{0}).
\end{align*}

\end{dem}

At this step, the objective of the next section will be to remove the
non-decreasingness condition used in the previous results.

\section{Extended real-valued minimax\ theorem\label{minimax1}}

We give in this section a minimax theorem for extended real-valued bifunctions
defined on the Cartesian product $X\times Y,$ for lcs $X$ and $Y.$

We will need some\ technical lemmas related to the continuity of convex
marginal functions.

\begin{lem}
\label{lemB1} Given\ a lsc convex function $F:X\times Y\rightarrow
\overline{\mathbb{R}}$ and a nonempty compact convex set $A\subset X,$ the
function\ $g:Y\rightarrow\mathbb{R}_{\infty}$ defined by\
\[
g(y):=\inf_{x\in A}F(x,y)
\]
is\ convex and lsc.
\end{lem}

\begin{dem}
The convexity of$\ g$\ is well-known as it is the marginal of the convex
function $F+\mathrm{I}_{A\times Y}$ (see, e.g., \cite[Tehorem 2.1.3$(v)$%
]{Za02}). To check that it is also lsc, we fix $\bar{y}\in Y$ and take an
arbitrary net $(y_{i})_{i}\subset Y$ that converges\ to $\bar{y}.$ We may
assume that $(y_{i})_{i}\subset\operatorname*{dom}g;$ otherwise,
$g(y_{i})=+\infty$ frequently for $i\in I,$ and\ the inequality $g(\bar
{y})\leq\liminf_{i}g(y_{i})$ obviously holds. We also take nets\ $(\alpha
_{i})_{i\in I},$ $(k_{i})_{i\in I}\subset\mathbb{R}_{+}$ such that $\alpha
_{i}\downarrow0$ and $k_{i}\uparrow+\infty.$ Then we find another net
$(x_{i})_{i\in I}\subset A$ such that
\begin{equation}
F(x_{i},y_{i})\leq\max\{g(y_{i}),-k_{i}\}+\alpha_{i},\text{ for all }i\in I.
\label{iaeb}%
\end{equation}
Moreover, taking into account the compactness of $A,$ we may assume that
$(x_{i})_{i}$ converges to some $\bar{x}\in A$ (w.o.l.g.~on $I$). Thus, taking
limits in (\ref{iaeb}),\ the lower semicontinuity\ of\ $F$\ gives rise to
\begin{align*}
F(\bar{x},\bar{y})  &  \leq\liminf_{i}F(x_{i},y_{i})\leq\liminf_{i}%
(\max\{g(y_{i}),-k_{i}\}+\alpha_{i})\\
&  =\liminf_{i}\max\{g(y_{i}),-k_{i}\}=\liminf_{i}g(y_{i}).
\end{align*}
Hence, since\ $\bar{x}\in A,$ we obtain\ $g(\bar{y})\leq F(\bar{x},\bar
{y})\leq\liminf_{i}g(y_{i}^{\ast})$ and the lower semicontinuity of $g$ at
$\bar{y}$ follows.
\end{dem}

We apply the previous lemma to a special case that interests us.

\begin{lem}
\label{lemB} Consider\ a function $f:X\times Y\rightarrow\overline{\mathbb{R}%
}$ and a nonempty compact convex set $A\subset X.$ If $f(\cdot,y)$ is concave
and usc, for each\ $y\in Y,$\ then the function\ $g:Y^{\ast}\rightarrow
\overline{\mathbb{R}}$ defined by\
\begin{equation}
g(y^{\ast}):=\inf_{x\in A}\left(  f(x,\cdot)\right)  ^{\ast}(y^{\ast})
\label{funcg}%
\end{equation}
is\ convex and lsc.
\end{lem}

\begin{dem}
Firstly, for each $y\in Y,$ the function $(x,y^{\ast})\in X\times Y^{\ast
}\mapsto\left\langle y,y^{\ast}\right\rangle -f(x,y)$ is convex and lsc, and
so is the pointwise supremum $F(x,y^{\ast}):=\sup_{y\in Y}\{\left\langle
y,y^{\ast}\right\rangle -f(x,y)\}=\left(  f(x,\cdot)\right)  ^{\ast}(y^{\ast
}).$ The conclusion\ follows then from Lemma \ref{lemB1}.\ 
\end{dem}

We give the main minimax result of this section for extended real-valued functions.

\begin{theo}
\label{MM1} Given a function $f:X\times Y\rightarrow\overline{\mathbb{R}}$ and
a nonempty convex set $A\subset X,$ we suppose\ the following
conditions\emph{:}

$(i)$ The set $A$ is compact.

$(ii)$ The functions\ $f(\cdot,y),$ $y\in Y,$ are concave and usc.$\ $\newline
Then we have
\begin{equation}
\inf_{y\in Y}\sup_{x\in A_{0}}f(x,y)\leq\max_{x\in A}\inf_{y\in Y}f(x,y),
\label{esd}%
\end{equation}
where
\begin{equation}
A_{0}:=\{x\in A:f(x,\cdot)\in\Gamma_{0}(Y)\}. \label{sao}%
\end{equation}

\end{theo}

\begin{dem}
First, let us observe that the function $\inf_{y\in Y}f(\cdot,y)$ is concave
and usc, so that $\sup_{x\in A}\inf_{y\in Y}f(x,y)=\max_{x\in A}\inf_{y\in
Y}f(x,y)$; that is, the supremum is attained. The inequality in (\ref{esd})
obviously holds whenever\ $A_{0}=\emptyset$, because in that case we would
have $\sup_{x\in A_{0}}f(x,y)=-\infty$ for all $y\in Y.$ Let us also
check\ that (\ref{esd}) holds when $A_{0}\not =\emptyset$ and
\[
\sup_{x\in A_{0}}f(x,y)=+\infty\ \text{for all }y\in Y.
\]
Indeed, in such a case, since $A_{0}\subset A$ we would have that $\sup_{x\in
A}f(x,y)=+\infty\ $for all $y\in Y.$ Thus, for each given $y\in Y,$ conditions
$(i)$-$(ii)$ yield an element $x(y)\in A$ such that $f(x(y),y)=+\infty.$ In
particular, the usc concave function $f(\cdot,y)$ is not proper; thus, it only
takes the infinite values$\ +\infty$ and $-\infty.$ But, from the definition
of the set $A_{0}$ (assumed nonempty), all\ the function $f(x,\cdot)$ for
$x\in A_{0}$ are proper, and so $f(x,y)=+\infty$ for all $x\in A_{0}$ and
$y\in Y.$ Consequently,
\[
\max_{x\in A}\inf_{y\in Y}f(x,y)=\sup_{x\in A}\inf_{y\in Y}f(x,y)\geq
\sup_{x\in A_{0}}\inf_{y\in Y}f(x,y)=+\infty=\inf_{y\in Y}\sup_{x\in A_{0}%
}f(x,y),
\]
and (\ref{esd}) holds in the current case too.

On account of the comments above, it suffices to prove (\ref{esd}) under the
assumption\ that $A_{0}\not =\emptyset$ and there exists\ some $y_{0}\in Y$
satisfying%
\begin{equation}
\sup_{x\in A_{0}}f(x,y_{0})<+\infty. \label{iii}%
\end{equation}
Then, using the definition of the conjugate,\textbf{ }we write
\begin{equation}
\max_{x\in A}\inf_{y\in Y}f(x,y)=\max_{x\in A}(-(f(x,\cdot))^{\ast}%
(\theta))=-\min_{x\in A}\left(  f(x,\cdot)\right)  ^{\ast}(\theta),
\label{tk2}%
\end{equation}
and Lemma \ref{lemB} together with\ the fact that $A_{0}\subset A$ yields
\[
\max_{x\in A}\inf_{y\in Y}f(x,y)=-\overline{\operatorname*{co}}\left(
\inf_{x\in A}\left(  f(x,\cdot)\right)  ^{\ast}\right)  (\theta)\geq
-\overline{\operatorname*{co}}\left(  \inf_{x\in A_{0}}\left(  f(x,\cdot
)\right)  ^{\ast}\right)  (\theta).
\]
Therefore, since\ $f(x,\cdot)\in\Gamma_{0}(Y)$ for all $x\in A_{0},$ and
$\sup_{x\in A_{0}}f(x,\cdot)\in\Gamma_{0}(Y)$ by (\ref{iii}),\ relation
(\ref{infsupcon}) implies\ that\
\[
\max_{x\in A}\inf_{y\in B}f(x,y)\geq-\left(  \sup_{x\in A_{0}}f(x,\cdot
)\right)  ^{\ast}(\theta)=\inf_{y\in Y}\sup_{x\in A_{0}}f(x,y),
\]
which is the desired inequality.
\end{dem}

If, in addition, the set $A$ defined in (\ref{sao})\ is such\ that
\begin{equation}
\inf_{y\in Y}\sup_{x\in A}f(x,y)\geq\inf_{y\in Y}\max_{x\in A}f(x,y),
\label{hyp0}%
\end{equation}
then Theorem \ref{MM1} implies the minimax equality%
\[
\max_{x\in A}\inf_{y\in Y}f(x,y)=\inf_{y\in Y}\max_{x\in A}f(x,y).
\]
In this line we have the following example where the convexity of $f(x,\cdot)$
is only required to hold for points\ $x\in\operatorname*{int}(A).$ A
finite-dimensional version of this result is also possible, up to replacing
the interior\ with the relative interior.$\ $

\begin{exam}
Given a function $f:X\times Y\rightarrow\overline{\mathbb{R}}$ and a nonempty
convex set $A\subset X,$ we assume the following conditions:

$(i)$ The set\ $A$ is compact.

$(ii)$ The functions\ $f(\cdot,y),$ $y\in Y,$ are concave and usc.

$(iii)$ $\left(  \operatorname*{int}(A)\right)  \cap\left\{  x\in
X:f(x,y)>-\infty\right\}  \neq\emptyset$ for all $y\in Y.$

$(iv)$ $\ f(x,\cdot)\in\Gamma_{0}(Y)$ for all $x\in\operatorname*{int}%
(A).$\newline Then we have the minimax equality :%
\[
\inf_{y\in Y}\max_{x\in A}f(x,y)=\max_{x\in A}\inf_{y\in Y}f(x,y).
\]
Indeed, on the one hand, condition $(iii)$ and the accessibility lemma (e.g.,
\cite[Lemma 1]{HLZ08}) ensure that
\[
\sup_{x\in\operatorname*{int}A}f(x,y)=\sup_{x\in A}f(x,y)\text{ for all }y\in
Y;
\]
that is,
\[
\inf_{y\in Y}\max_{x\in A}f(x,y)=\inf_{y\in Y}\sup_{x\in\operatorname*{int}%
A}f(x,y).
\]
But $(iv)$ implies that $\operatorname*{int}(A)\subset A_{0}=\{x\in
A:f(x,\cdot)\in\Gamma_{0}(Y)\}\subset A,$ and so Theorem \ref{MM1} yields
\[
\inf_{y\in Y}\max_{x\in A}f(x,y)=\inf_{y\in Y}\sup_{x\in\operatorname*{int}%
A}f(x,y)\leq\max_{x\in A}\inf_{y\in Y}f(x,y)\leq\inf_{y\in Y}\max_{x\in
A}f(x,y),
\]
proving the desired minimax equality.
\end{exam}

A\ localized version of Theorem \ref{MM1} is given in the corollary below. Its
proof is immediate by applying this theorem to the new bifunction%
\[
(x,y)\longmapsto\tilde{f}(x,y):=f(x,y)+\mathrm{I}_{B}(y),
\]
which is easily shown to satisfy\ conditions $(i)$-$(ii)$ of Theorem
\ref{MM1}, and to the set
\begin{equation}
\tilde{A}_{0}:=\{x\in A:\tilde{f}(x,\cdot)\in\Gamma_{0}(Y)\}=\{x\in
A:f(x,\cdot)+\mathrm{I}_{B}(\cdot)\in\Gamma_{0}(Y)\}. \label{saob}%
\end{equation}
Note that here the convexity of the given set $B$ is implicit, because the
nonemptyness of the set $\tilde{A}_{0}$ guarantees the convexity
of\ $B\cap\operatorname*{dom}f(x,\cdot)$ for all $x\in\tilde{A}_{0}.$

\begin{cor}
\label{MM} Given a function $f:X\times Y\rightarrow\overline{\mathbb{R}}$ and
nonempty\ sets $A\subset X,$ $B\subset Y,$ we assume the following conditions:

$(i)$ The set\ $A$ is convex and compact.

$(ii)$ The functions\ $f(\cdot,y),$ $y\in B,$ are concave and usc.\ \newline
Then we have
\begin{equation}
\inf_{y\in B}\sup_{x\in A_{0}}f(x,y)\leq\max_{x\in A}\inf_{y\in B}f(x,y),
\label{esdb}%
\end{equation}
where $A_{0}:=\{x\in A:f(x,\cdot)+\mathrm{I}_{B}(\cdot)\in\Gamma_{0}(Y)\}.$
\end{cor}

\section{Minimax theorem\label{minimax2}}

We give another variant of the minimax theorem, dropping out\ the lower
semicontinuity condition of the functions $f(x,\cdot),$ $x\in A,$ used in
Theorem \ref{MM1}. Instead, we use here the condition that the function $f$ is
finite-valued on the set $A\times B.$ As in the previous section, we also
assume here that $X$ and $Y$ are two lcs.

\begin{theo}
\label{MMB}\ Given a function $f:X\times Y\rightarrow\overline{\mathbb{R}}$
and nonempty convex sets $A\subset X,$ $B\subset Y$ such that $A\times
B\subset f^{-1}(\mathbb{R}),$ we assume the following conditions\emph{:}

$(i)$ The set\ $A$ is compact.

$(ii)$ The functions\ $f(\cdot,y)$, $y\in B,$ are concave and usc.\ \newline
Then, we have%
\[
\inf_{y\in B}\sup_{x\in A_{1}}f(x,y)\leq\max_{x\in A}\inf_{y\in B}f(x,y),
\]
where%
\[
A_{1}:=\{x\in A:f(x,\cdot)\text{\ is convex}\}.
\]

\end{theo}

\begin{dem}
First, note that the relation $A\times B\subset f^{-1}(\mathbb{R})$ together
with condition $(ii)$ entails%
\[
\sup_{x\in A}f(x,y)=\max_{x\in A}f(x,y)<+\infty\text{ \ for every }y\in B,
\]
implying\ that%
\begin{equation}
B\subset\operatorname*{dom}\left(  \sup_{x\in A}f(x,\cdot)\right)  .
\label{hypp}%
\end{equation}
Next, we introduce the family%
\[
\mathcal{F}^{B}:=\{L\subset Y:L\text{ is a finite-dimensional linear subspace
that intersects }B\},
\]
and pick an $L\in\mathcal{F}^{B}.\ $Arguing as in (\ref{tk2}) and using
(\ref{tricy}),\ we\ write\
\begin{align}
\max_{x\in A}\inf_{y\in L\cap B}f(x,y)  &  =\max_{x\in A}\inf_{y\in
Y}(f(x,y)+\mathrm{I}_{L\cap B}(y))\nonumber\\
&  =\max_{x\in A}\left[  -(f(x,\cdot)+\mathrm{I}_{L\cap B})^{\ast}%
(\theta)\right] \nonumber\\
&  =\max_{x\in A}\left[  -\left(  \operatorname*{cl}\nolimits_{y}%
(f(x,\cdot)+\mathrm{I}_{L\cap B}(\cdot))\right)  ^{\ast}(\theta)\right]
\nonumber\\
&  =-\inf_{x\in A}\left[  \left(  \operatorname*{cl}\nolimits_{y}%
(f(x,\cdot)+\mathrm{I}_{L\cap B}(\cdot))\right)  ^{\ast}(\theta)\right]  ,
\label{tk2b}%
\end{align}
where $\operatorname*{cl}\nolimits_{y}(f(x,\cdot)+\mathrm{I}_{L\cap B}%
(\cdot))$ denotes the lsc hull\ of the function $f(x,\cdot)+\mathrm{I}_{L\cap
B}(\cdot)$ with respect\ to the variable $y.$ Furthermore,\ since the function
$g_{L}:Y^{\ast}\rightarrow\overline{\mathbb{R}},$ defined by
\[
g_{L}(y^{\ast}):=\inf_{x\in A}\left(  \operatorname*{cl}\nolimits_{y}%
(f(x,\cdot)+\mathrm{I}_{L\cap B}(\cdot))\right)  ^{\ast}(y^{\ast})=\inf_{x\in
A}\left(  f(x,\cdot)+\mathrm{I}_{L\cap B}(\cdot)\right)  ^{\ast}(y^{\ast}),
\]
is convex and lsc by Lemma \ref{lemB}, and\ $A_{1}\subset A$ obviously, the
inequality\ in (\ref{tk2b}) reads\
\begin{align}
\max_{x\in A}\inf_{y\in L\cap B}f(x,y)  &  =-\overline{\operatorname*{co}%
}\left(  \inf_{x\in A}\left(  \operatorname*{cl}\nolimits_{y}(f(x,\cdot
)+\mathrm{I}_{L\cap B}(\cdot))\right)  ^{\ast}\right)  (\theta)\nonumber\\
&  \geq-\overline{\operatorname*{co}}\left(  \inf_{x\in A_{1}}\left(
\operatorname*{cl}\nolimits_{y}(f(x,\cdot)+\mathrm{I}_{L\cap B}(\cdot
))\right)  ^{\ast}\right)  (\theta). \label{yah}%
\end{align}
Note\ that the convex function\ $\operatorname*{cl}\nolimits_{y}%
(f(x,\cdot)+\mathrm{I}_{L\cap B}(\cdot))$ above is (lsc and) proper
because\ the convex function $f(x,\cdot)+\mathrm{I}_{L\cap B}(\cdot)$ is
proper, thanks to (\ref{hypp}), and has a finite-dimensional effective domain
(see \cite{BV10}). As a consequence of that, the function
\[
\varphi:=\sup_{x\in A_{1}}\left(  \operatorname*{cl}\nolimits_{y}%
(f(x,\cdot)+\mathrm{I}_{L\cap B}(\cdot))\right)
\]
does not take the value $-\infty.$ Moreover, by (\ref{hypp}) we have\
\begin{equation}
\emptyset\neq L\cap B\subset\operatorname*{dom}\left(  \max_{x\in A}%
f(x,\cdot)+\mathrm{I}_{L\cap B}(\cdot)\right)  \subset\operatorname*{dom}%
\left(  \sup_{x\in A_{1}}\left(  f(x,\cdot)+\mathrm{I}_{L\cap B}%
(\cdot)\right)  \right)  , \label{yh}%
\end{equation}
which shows that
\begin{equation}
\operatorname*{dom}\left(  \sup_{x\in A_{1}}\left(  f(x,\cdot)+\mathrm{I}%
_{L\cap B}(\cdot)\right)  \right)  =\operatorname*{dom}\left(  f(x,\cdot
)+\mathrm{I}_{L\cap B}(\cdot)\right)  =L\cap B\text{ \ for all }x\in A_{1}.
\label{yh2}%
\end{equation}
In other words, since $\varphi\leq\sup_{x\in A_{1}}\left(  f(x,\cdot
)+\mathrm{I}_{L\cap B}(\cdot)\right)  $, we have $L\cap B\subset
\operatorname*{dom}\varphi$ and the function $\varphi$ is proper.
Consequently, (\ref{infsupcon}) entails\
\[
\varphi^{\ast}(\theta)=\overline{\operatorname*{co}}\left(  \inf_{x\in A_{1}%
}\left(  \operatorname*{cl}\nolimits_{y}(f(x,\cdot)+\mathrm{I}_{L\cap B}%
(\cdot))\right)  ^{\ast}\right)  (\theta),
\]
and (\ref{yah}) gives rise, according to Moreau's theorem, to
\begin{equation}
\max_{x\in A}\inf_{y\in L\cap B}f(x,y)\geq-\left(  \sup_{x\in A_{1}}\left(
\operatorname*{cl}\nolimits_{y}(f(x,\cdot)+\mathrm{I}_{L\cap B}(\cdot
))\right)  \right)  ^{\ast}(\theta)=\inf_{y\in Y}\varphi(y). \label{tog}%
\end{equation}
Moreover, due to (\ref{yh2}), for each $x\in A_{1}$ the set
$\operatorname*{dom}(f(x,\cdot)+\mathrm{I}_{L\cap B}(\cdot))=L\cap B$ is
finite-dimensional and, therefore,\ (\ref{amin}) entails\
\[
\varphi=\sup_{x\in A_{1}}\left(  \operatorname*{cl}\nolimits_{y}%
(f(x,\cdot)+\mathrm{I}_{L\cap B}(\cdot))\right)  =\operatorname*{cl}%
\nolimits_{y}\left(  \sup_{x\in A_{1}}\left(  f(x,\cdot)+\mathrm{I}_{L\cap
B}(\cdot)\right)  \right)  .
\]
In other words, (\ref{tog}) and (\ref{ehuls}) yield\
\begin{align}
\max_{x\in A}\inf_{y\in L\cap B}f(x,y)  &  \geq\inf_{y\in Y}\left(
\operatorname*{cl}\nolimits_{y}\left(  \sup_{x\in A_{1}}\left(  f(x,\cdot
)+\mathrm{I}_{L\cap B}(\cdot)\right)  \right)  \right) \nonumber\\
&  =\inf_{y\in Y}\sup_{x\in A_{1}}\left(  f(x,\cdot)+\mathrm{I}_{L\cap
B}(\cdot)\right) \nonumber\\
&  \geq\inf_{y\in Y}\sup_{x\in A_{1}}\left(  f(x,\cdot)+\mathrm{I}_{B}%
(\cdot)\right)  =\inf_{y\in B}\sup_{x\in A_{1}}f(x,\cdot). \label{sab}%
\end{align}
Let us, finally, endow the family $\mathcal{F}^{B}$ with the partial order
given by ascending inclusions, \
\[
L_{1}\preccurlyeq L_{2},\text{ }L_{1},L_{2}\in\mathcal{F}^{B}\Leftrightarrow
L_{1}\subset L_{2}.
\]
Then, applying Corollary\ \ref{mmmonotone} in $(\mathcal{F}^{B},\preccurlyeq)$
to\ the non-increasing net of the usc (concave) functions $\varphi_{L}%
:=\inf_{y\in L\cap B}f(\cdot,y),$ $L\in\mathcal{F}^{B},$ (\ref{sab}) yields
\begin{align*}
\max_{x\in A}\inf_{y\in B}f(x,y)  &  =\max_{x\in A}\inf_{L\in\mathcal{F}^{B}%
}\inf_{y\in L\cap B}f(\cdot,y)\\
&  =\max_{x\in A}\inf_{L\in\mathcal{F}^{B}}\varphi_{L}(x)=\inf_{L\in
\mathcal{F}^{B}}\max_{x\in A}\varphi_{L}(x)\\
&  \geq\inf_{L\in\mathcal{F}^{B}}\inf_{y\in B}\sup_{x\in A_{1}}f(x,\cdot
)=\inf_{y\in B}\sup_{x\in A_{1}}f(x,\cdot),
\end{align*}
and we are done with the proof.
\end{dem}

The classical\ minimax theorem straightforwardly\ follows from Theorem
\ref{MMB} when all the functions $f(x,\cdot),$ $x\in A,$ are convex (that is,
when $A_{1}=A$). Next, we give a useful application of Theorems \ref{MM1} and
\ref{MMB}. Given $n\geq1,$ we denote
\[
\Delta_{n}:=\{\lambda:=(\lambda_{1},\ldots,\lambda_{n})\in\mathbb{R}%
^{n}:\lambda_{k}\geq0,\text{ }%
{\textstyle\sum\nolimits_{1\leq k\leq n}}
\lambda_{k}=1\}.
\]

\begin{cor}
\label{minimax} Given\ a collection of convex functions $f_{k}:X\rightarrow
\overline{\mathbb{R}},$ $1\leq k\leq n,$ we write $f:=\max_{1\leq k\leq
n}f_{k}$ and suppose that $\operatorname*{dom}f\neq\emptyset.$ We assume that
at least one of the following conditions holds:

$(i)$ All the $f_{k}$'s are proper.

$(ii)$ All the $f_{k}$'s are lsc.\newline Then, we have that
\[
\inf_{x\in X}f(x)=\max_{\lambda\in\Delta_{n}}\inf_{x\in X}\text{ }%
{\textstyle\sum\limits_{1\leq k\leq n}}
\lambda_{k}f_{k}(x).
\]

\end{cor}

\begin{dem}
We consider the function $F:\mathbb{R}^{n}\times X\rightarrow\overline
{\mathbb{R}}$ defined as
\begin{equation}
F(\lambda,x):=f_{\lambda}(x)-\mathrm{I}_{\mathbb{R}_{+}^{n}}(\lambda),\text{
}\lambda\in\mathbb{R}^{n},\text{ }x\in X, \label{phi}%
\end{equation}
where $f_{\lambda}:=%
{\textstyle\sum\nolimits_{1\leq k\leq n}}
\lambda_{k}f_{k}.$ We also\ denote $A:=\Delta_{n}\subset\mathbb{R}^{n},$
$B:=\operatorname*{dom}f\subset X,$ and
\[
A_{0}:=\{\lambda\in\Delta_{n}:\lambda_{k}>0\text{ for all }k=1,\ldots,n\}.
\]
Hence, the functions $F(\cdot,x),$ $x\in B,$\ are usc and concave (indeed,
affine), whereas the functions $F(\lambda,\cdot),$ $\lambda\in A,$ are
convex.\ Also, it is clear that the set $A$ is convex and compact.

Let us first assume that all the $f_{k}$'s are proper, so that $A\times
B\subset F^{-1}(\mathbb{R})$. Thus, Theorem \ref{MMB}\ applies and gives us\
\[
\max_{\lambda\in A}\inf_{x\in B}F(\lambda,x)=\inf_{x\in B}\max_{\lambda\in
A}F(\lambda,x),
\]
showing that
\[
\max_{\lambda\in\Delta_{n}}\inf_{x\in X}f_{\lambda}(x)=\max_{\lambda\in
\Delta_{n}}\inf_{x\in\operatorname*{dom}f}F(\lambda,x)=\inf_{x\in
\operatorname*{dom}f}\max_{\lambda\in\Delta_{n}}F(\lambda,x)=\inf_{x\in
X}f(x),
\]
where the last equality comes from the definition of the maximum function.

Second, if all the $f_{k}$'s are lsc, then the functions $F(\lambda
,\cdot)+\mathrm{I}_{B}$ ($=F(\lambda,\cdot),$ as $B=\operatorname*{dom}%
F(\lambda,\cdot)$), $\lambda\in A_{0},$ belong to $\Gamma_{0}(X).$\ Thus, by
Theorem \ref{MM1}, we obtain that
\[
\max_{\lambda\in A}\inf_{x\in B}F(\lambda,x)\geq\inf_{x\in B}\sup_{\lambda\in
A_{0}}F(\lambda,x).
\]
Moreover, we have that
\[
\inf_{x\in B}\sup_{\lambda\in A_{0}}F(\lambda,x)=\inf_{x\in\operatorname*{dom}%
f}\sup_{\lambda\in A_{0}}f_{\lambda}(x)=\inf_{x\in\operatorname*{dom}f}%
\sup_{\lambda\in\Delta_{n}}f_{\lambda}(x),
\]
and, as above, we deduce that%
\[
\inf_{x\in X}f(x)=\inf_{x\in\operatorname*{dom}f}\max_{\lambda\in\Delta_{n}%
}F(\lambda,x)\leq\max_{\lambda\in\Delta_{n}}\inf_{x\in\operatorname*{dom}%
f}F(\lambda,x)=\max_{\lambda\in\Delta_{n}}\inf_{x\in X}F(\lambda,x)\leq
\inf_{x\in X}f(x),
\]
which in turn leads us to the desired conclusion.
\end{dem}

Corollary \ref{minimax} easily allows us to formulate the subdifferential of
the maximum function $f:=\max_{1\leq k\leq n}f_{k},$ providing a new proof and
a slight extension\ of \cite[Corollary 2.8.11]{Za02} (see, also, references
therein to trace back the origin of this result) to improper functions.

\begin{cor}
\label{cor} With the assumptions of Corollary \ref{minimax}, for every $x\in
X$ and $\varepsilon\geq0$ we have that
\begin{equation}
\partial_{\varepsilon}f(x)=%
{\textstyle\bigcup\limits_{\lambda\in\Delta_{n}}}
\partial_{\left(  \varepsilon+f_{\lambda}(x)-f(x)\right)  }f_{\lambda}(x),
\label{exp}%
\end{equation}
where $f_{\lambda}:=%
{\textstyle\sum\nolimits_{1\leq k\leq n}}
\lambda_{k}f_{k}.$ In particular, for $\varepsilon=0$ we have
\begin{equation}
\partial f(x)=%
{\textstyle\bigcup}
\left\{  \partial f_{\lambda}(x):\lambda\in\Delta_{n},\text{ }f_{\lambda
}(x)=f(x)\right\}  . \label{expb}%
\end{equation}

\end{cor}

\begin{dem}
We\ fix $x_{0}\in X$\ and $\varepsilon\geq0.$ Formula (\ref{expb}) is an
immediate consequence of (\ref{exp}), due to the fact $\partial_{\left(
f_{\lambda}(x_{0})-f(x_{0})\right)  }f_{\lambda}(x_{0})$ is empty whenever
$f_{\lambda}(x_{0})<f(x_{0}).$ Thus, we only need to prove the inclusion
\textquotedblleft$\subset$\textquotedblright\ in (\ref{exp}) because the
opposite inclusion there can be\ easily checked. Let us first suppose that
$\theta\in\partial_{\varepsilon}f(x_{0})$ or, equivalently, that
\[
f(x_{0})\leq\inf_{x\in X}f(x)+\varepsilon.
\]
Observe that $f(x_{0})\in\mathbb{R}.$ Therefore, according to Corollary
\ref{minimax}, there exists some $\bar{\lambda}\in\Delta_{n}$ such that
\begin{align*}
f(x_{0})  &  \leq\inf_{x\in X}\max_{\lambda\in\Delta_{n}}f_{\lambda
}(x)+\varepsilon\\
&  =\max_{\lambda\in\Delta_{n}}\inf_{x\in X}f_{\lambda}(x)+\varepsilon
=\inf_{x\in X}f_{\bar{\lambda}}(x)+\varepsilon.
\end{align*}
In particular, we have that $f_{\bar{\lambda}}(x_{0})\leq f(x_{0})\leq
f_{\bar{\lambda}}(x_{0})+\varepsilon$ and, so,
\[
f_{\bar{\lambda}}(x_{0})=f(x_{0})+f_{\bar{\lambda}}(x_{0})-f(x_{0})\leq
\inf_{x\in X}f_{\bar{\lambda}}(x)+f_{\bar{\lambda}}(x_{0})-f(x_{0}%
)+\varepsilon;
\]
that is, $\theta\in\partial_{(f_{\bar{\lambda}}(x_{0})-f(x_{0})+\varepsilon
)}f_{\bar{\lambda}}(x_{0}).$

More generally, if $x^{\ast}\in\partial_{\varepsilon}f(x_{0}),$ then
$\theta\in\partial_{\varepsilon}(f-\left\langle x_{0}^{\ast},\cdot
\right\rangle )(x_{0})$ and we apply the paragraph above to the convex
functions $\tilde{f}_{k}:=f_{k}-\left\langle x_{0}^{\ast},\cdot\right\rangle
.$
\end{dem}

We close the paper with the following corollary to show that\ the Moreau
theorem (see (\ref{moreau})) can also be obtained from the minimax theorem,
Theorem \ref{MMB}. \ \ This proves that somehow these two\ results\ can be considered\ equivalent.

\begin{cor}
For every function $f\in\Gamma_{0}(X),$ we have that $f^{\ast\ast}=f.$
\end{cor}

\begin{dem}
Given a function\ $f\in\Gamma_{0}(X),$\ we fix $x\in X$ and a $\theta
$-neighborhood $U\subset X.$ Then, by definition of the biconjugate, we write
\[
f^{\ast\ast}(x)=\sup_{x^{\ast}\in X^{\ast}}\inf_{y\in\operatorname*{dom}%
f}g_{x^{\ast}}(y)\geq\sup_{x^{\ast}\in U^{\circ}}\inf_{y\in\operatorname*{dom}%
f}g_{x^{\ast}}(y),
\]
where the functions $g_{x^{\ast}}\in\Gamma_{0}(X),$ $x^{\ast}\in X^{\ast},$
are defined by
\[
g_{x^{\ast}}(y):=\left\langle x^{\ast},x-y\right\rangle +f(y).
\]
Observe that the functions $x^{\ast}\mapsto g_{y}(x^{\ast}),$ $y\in X,$ are
concave (and usc). Also, the convex set $U^{\circ}$ is $w^{\ast}$-compact
thanks to Dieudonné's Theorem. Therefore, applying Theorem\ \ref{MM1}$\ $with
$A:=U^{\circ}$\ and $B:=X\ $gives rise to%
\[
f^{\ast\ast}(x)\geq\sup_{x^{\ast}\in U^{\circ}}\inf_{y\in X}g_{x^{\ast}%
}(y)=\inf_{y\in X}\sup_{x^{\ast}\in U^{\circ}}\{\left\langle x^{\ast
},x-y\right\rangle +f(y)\},
\]
which in turn yields, using\ Example \ref{exam},
\[
f^{\ast\ast}(x)\geq\sup_{U\in\mathcal{N}_{X}}\inf_{y\in X}\sup_{x^{\ast}\in
U^{\circ}}\{\left\langle x^{\ast},x-y\right\rangle +f(y)\}=f(x).
\]
The proof is finished because\ the inequality $f^{\ast\ast}(x)\leq f(x)$
always holds.
\end{dem}

\end{document}